\documentclass[conference]{IEEEtran}

\makeatletter
\def\ps@headings{%
\def\@oddhead{\mbox{}\scriptsize\rightmark \hfil \thepage}%
\def\@evenhead{\scriptsize\thepage \hfil \leftmark\mbox{}}%
\def\@oddfoot{}%
\def\@evenfoot{}}
\makeatother
\usepackage{graphicx}
\usepackage[center]{caption}
\usepackage{epsfig,latexsym,amsmath,amsfonts}
\usepackage{amssymb}
\usepackage{placeins}
\usepackage{subfigure}
\usepackage{verbatim}
\usepackage{cite,url}
\usepackage{epsf,bm}

\begin{document}

\title{A Soft Sensing-Based Cognitive Access Scheme Exploiting Primary Feedback}

\author{\large Ahmed M. Arafa$^\dagger$, Karim G. Seddik$^\ddagger$, Ahmed K. Sultan$^*$, Tamer ElBatt$^\dagger$ and Amr A. El-Sherif$^*$  \\ [.1in]
\small  \begin{tabular}{c} $^\dagger$Wireless Intelligent Networks Center (WINC), Nile University, Smart Village, Egypt.\\
$^\ddagger$Electronics Engineering Department, American University in Cairo, AUC Avenue, New Cairo 11835, Egypt.\\
$^*$Department of Electrical Engineering, Alexandria University, Alexandria 21544, Egypt. \\
email: ahmed.arafa@nileu.edu.eg, kseddik@aucegypt.edu, salatino@stanfordalumni.edu,  telbatt@ieee.org, amr.elsherif@ieee.org
\end{tabular} }

\maketitle

\begin{abstract}

In this paper, we examine a cognitive spectrum access scheme in which secondary users exploit the primary feedback information. We consider an overlay secondary network employing a random access scheme in which secondary users access the channel by certain access probabilities that are function of the spectrum sensing metric. In setting our problem, we assume that secondary users can eavesdrop on the primary link's feedback. We study the cognitive radio network from a queuing theory point of view. Access probabilities are determined by solving a secondary throughput maximization problem subject to a constraint on the primary queues' stability. First, we formulate our problem which is found to be non-convex. Yet, we solve it efficiently by exploiting the structure of the secondary throughput equation. Our scheme yields improved results in, both, the secondary user throughput and the primary user packet delay. In addition, it comes very close to the \textit{optimal} genie-aided scheme in which secondary users act upon the presumed perfect knowledge of the primary user's activity\footnotemark.

\end{abstract}

\footnotetext{This work was supported by a grant from the Egyptian National Telecommunications Regulatory Authority (NTRA).}
\footnotetext[2]{Tamer ElBatt is also affiliated with the EECE Dept., Faculty of Engineering, Cairo University.}

\section{Introduction}\label{Int}

Cognitive radio technology is a communication paradigm that emerged in order to solve the spectrum scarcity problem by allowing unlicensed (or secondary) users to exploit the under-utilized spectrum of the licensed (or primary) users. Coexistence of such secondary users along with primary ones is allowed provided that minimal or no harm is caused upon the primary network, and that a minimum quality of service is guaranteed for primary users.

In a typical cognitive radio setting, the cognitive transmitter senses the primary activity and decides on accessing the channel based on the sensing outcome. This setting is problematic in the sense that cognitive users are not aware of their impact on the primary network, besides the usual sensing errors. This, in turn, induced two ideas to alleviate these hurdles. The first is to decrease the sensing error rate, as in \cite{Soft_Sensing}, where the authors introduced a novel design in which the value of the test statistic is used as a confidence measure for the sensing outcome. This value is then used to specify a channel access probability for the secondary network. The access probabilities as a function of the sensing metric are obtained by solving an optimization problem formulated to maximize the secondary throughput given a constraint on the primary queue stability. This idea is known as \textit{soft sensing}, and was first introduced in \cite{jafar}, however, the focus was on physical layer power adaptation to maximize the capacity of the secondary link.

The second idea proposed in the literature is to make the secondary user aware of the primary activity by leveraging the feedback sent from the primary receiver to the primary transmitter and optimizing its transmission strategy based on its effect on the primary receiver. For instance, in \cite{bits}, the secondary user observes the automatic repeat request (ARQ) from the primary receiver. The ARQs reflect the primary user's achieved packet rate. The cognitive radio's objective is to maximize the secondary throughput under the constraint of guaranteeing a certain packet rate for the primary user. In \cite{fabio}, the authors use a partially observable Markov decision process (POMDP) to devise an optimized admission control policy. Secondary power control on the basis of the primary link feedback is investigated in \cite{power_control_primary_arq}. In \cite{marco2}, the optimal transmission policy for the secondary user, when the primary user adopts a retransmission based error control scheme, is investigated. The policy of the secondary user determines how often it transmits according to the retransmission state of the packet being served by the primary user. The resulting optimal strategy of the secondary user is proven to have a unique structure. In particular, the optimal throughput is achieved by the secondary user by concentrating its interference to the primary user in the first transmission attempt of a packet. A simple idea is introduced in a previous work \cite{Feedback_SPAWC} in which secondary users refrain from accessing the channel upon hearing a NACK from the primary receiver allowing for an interference-free primary retransmission, thereby increasing secondary throughput and decreasing primary packet delay.

In this paper, we introduce a hybrid scheme in which we capture the benefits of a feedback-based access scheme introduced on top of soft sensing, and accordingly, high sensing reliability is attained besides awareness of the primary environment. We consider a secondary network employing a random access scheme in which secondary users access the channel by certain access probabilities that are function of the sensing metric. The network is studied from a queuing theoretic perspective, and access probabilities are determined by solving an optimization problem subject to a constraint on the primary user's queue stability. In addition, secondary users can overhear and, hence, leverage the primary link's feedback; secondary users back-off completely from accessing the channel upon hearing a NACK, and attempt accessing if an ACK/no feedback is overheard. This leads to significant improvements in the secondary user's throughput as well as the primary user's packet delay. This is attributed to the high sensing reliability, due to the use of soft sensing, as well as avoiding sure collisions between primary and secondary users when secondary users back-off upon hearing a NACK. Our scheme is shown to outperform both soft sensing and conventional hard decision sensing not leveraging feedback information and approaches the \textit{optimal} genie-aided scheme in which secondary users have perfect knowledge of the primary users' activity and, hence, make best use of it.

The rest of the paper is organized as follows. The system model is presented in Section \ref{sysmod}. A background on the soft sensing scheme is presented in Section \ref{SoftSense}. The proposed feedback- and soft sensing-based access scheme is described and analyzed in Section \ref{FbScheme}. Performance results are given in Section \ref{Simulations}. Finally, conclusions are drawn in Section \ref{Concl}.

\section{System Model}\label{sysmod}

We consider the uplink of a TDMA system consisting of $M_p$ primary users (PU), along which we have $M_s$ secondary users (SU) attempting to access the channel using a Slotted ALOHA scheme. Let $\mathcal{M}_p=\{1,2,...,M_p\}$ denote the set of all primary users, and $\mathcal{M}_s=\{1,2,...,M_s\}$ denote the set of all secondary users.

We consider an overlay system in which secondary users attempt to send their packets only when primary users are sensed to be idle. We adopt a collision model for interference whereby packets are lost if more than one transmission proceed at a time. At the beginning of each time slot, each SU senses the channel and if found idle, it accesses the channel with a certain access probability. Simple energy detection \cite{akyildiz2006next} is adopted as the sensing mechanism since it does not need prior information of the PU signal or its structure.

The channel is modeled as a Rayleigh flat fading channel with additive white Gaussian noise (AWGN.) Thus, the received signal at node $j$ from node $q$ at time slot $t$ is given by
\begin{equation}\label{eqn: Rx Signal}
y_{qj}^t=\sqrt{G_qr_{qj}^{-\gamma}}h_{qj}^tx_q^t+n_j^t,
\end{equation}
where $G_q$ is the transmitted power, $r_{qj}$ is the distance between the two nodes, and $\gamma$ is the path loss exponent. $x^t_q$ is the transmitted signal, which is assumed to be drawn from any constant modulus constellation, M-ary PSK for instance, with zero mean and unit variance. $h^t_{qj}$ is the channel coefficient between the two nodes, modeled as i.i.d. circularly symmetric complex Gaussian random variable with zero mean and unit variance. The noise term $n^t_j$ is also modeled as i.i.d. circularly symmetric complex Gaussian random variable with zero mean and variance $N_0$. We assume the channel is stationary and independent from slot to slot, thus, the superscript $t$ is dropped in the rest of this paper.

For a transmission to be successful, the channel must not be in outage, i.e. the received SNR should not be smaller than a pre-specified threshold $\zeta$. From the signal model in (\ref{eqn: Rx Signal}), the outage probability between nodes $q$ and $j$ is given by $P_{qj}^o=Pr\left\{|h_{qj}|^2<\frac{\zeta N_0 r_{qj}^\gamma}{G_i}\right\}=1-\exp\left(-\frac{\zeta N_0 r_{qj}^\gamma}{G_i}\right)$.


Each primary user has an infinite buffer for storing its incoming packets. The packet arrival processes at primary queues are assumed to be Bernoulli i.i.d. with an average arrival rate of $\lambda_q$ for user $q$. A slot duration is equal to the packet transmission time, and therefore, $0\leq\lambda_q\leq1$, $\forall q$. Assuming symmetry, for mathematical tractability, all $\lambda_q$'s are the same for all primary users and are equal to $\lambda_p$. Furthermore, in our model, we consider the case where secondary users always have packets to send.

Primary users access the channel by dividing the channel resources, time in this case, among them; hence, each node is allocated a fraction of the time. Let $\Omega_p=[\omega_p^1,\omega_p^2,...,\omega_p^{M_p}]$ denote a resource-sharing vector, where $\omega_p^q\geq0$ is the fraction of time allocated to node $q\in \mathcal{M}_p$, or it can
represent the probability that node $i$ is allocated the time slot \cite{M.Kobayashi}. Therefore, the set of all feasible resource sharing vectors is specified as $\digamma_p=\left\{ \Omega_p=(\omega_p^1,\omega_p^2,...,\omega_p^{M_p})\in \mathbb{R}^{+M_p}:\sum_{q\in M_p} \omega_p^q \leq 1 \right\}$, where $\mathbb{R}^{+M_p}$ is the set of $M_p$ dimensional vectors with real, non-negative elements.

In the proposed model, we leverage an error-free primary feedback channel via which the primary receiver sends a feedback by the end of each time slot to acknowledge the reception of packets. Accordingly, an ACK is sent if a packet is correctly received, and a NACK is sent if a packet is lost. Failure of reception is attributed to either primary channel outage, or collision between secondary and primary packets. In case of an idle slot, no feedback is sent. Secondary users are assumed to overhear this primary feedback perfectly and act as follows: if an ACK/no feedback is heard, the secondary users behave normally, and start sensing the channel in the next time slot. On the other hand, if a NACK is heard, all secondary users back-off in the next time slot allowing for an interference-free retransmission of the erroneous primary packet. Accordingly, sure collisions can be avoided since the reception of a NACK triggers the PU to send in the next time slot with probability one.

In the sequel, we assume symmetry conditions, for simplicity of analysis and presentation, in which all primary users' transmit powers are equal and all distances between secondary and primary users are equal. Therefore the subscript $qj$ is dropped in the rest of the paper.

In the next section, we present a background on the so-called soft sensing scheme which was briefly mentioned in Section \ref{Int}.

\section{Background: Soft Sensing-Based Access}\label{SoftSense}

\begin{figure}
\includegraphics[scale=0.3]{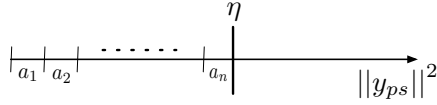}
\centering
\caption{Soft Sensing: Division of the interval $[0,\eta]$ into subintervals and their associated access probabilities.}\label{fig: access_intervals}
\vspace{-0.15in}
\end{figure}

We focus on the concept of soft sensing originally introduced in \cite{Soft_Sensing} which basically uses the energy statistic $\left\|y_{ps}\right\|^2$ acquired from the energy detector as a measure of reliability, where subscript $p$ denotes the primary user and $s$ denotes the secondary user. The lower the value of $\left\|y_{ps}\right\|^2$ compared to the decision threshold $\eta$, the more certain secondary users become that primary users are idle in the time slot in question. This observation is exploited to yield the powerful concept of soft sensing as follows:
\begin{itemize}
\item The interval $[0,\eta]$ is divided into $n$ subintervals as shown in Fig. \ref{fig: access_intervals}.
\item For each subinterval $i \in [1,n]$, an access probability $a_i$ is assigned.
\item If $\left\|y_{ps}\right\|^2$ lies in the $i^{th}$ subinterval, the SU attempts accessing the channel with probability $a_i$.
\item If $\left\|y_{ps}\right\|^2$ value is greater than $\eta$, the SU does not access the channel.
\end{itemize}

Intuition may suggest that access probabilities associated with subintervals far less than the threshold $\eta$ are given higher values than those associated with ones near the threshold, i.e. their values are sorted in a descending order as $i$ goes from $1$ to $n$. This is mainly because of the very low probability of collision with a PU packet whenever the energy statistic lies in the subintervals close to zero. On the other hand, there is a higher risk of collision whenever the energy statistic lies in the subintervals close to the energy threshold.

In this work, the stability of the PU queue is studied as the performance measure. Access probabilities are chosen such that the SU throughput is maximized provided that the PU queue is stable. Stability can be loosely defined as keeping a quantity of interest bounded, in this case, the queue size. For a more general and rigorous definition of stability, see \cite{Rao88} and \cite{Szpan94}. If the arrival and service processes of a queuing system are strictly stationary, one can apply Loynes' theorem to check for stability \cite{Loynes}. This theorem states that if the average arrival rate is less than the average service rate of a queuing system whose arrival and service processes are strictly stationary, then the queue is stable, otherwise it is unstable.

Therefore, the baseline problem, without feedback, can be formulated as maximizing the secondary throughput subject to the primary queue being stable. That is,
\begin{equation}\label{eqn: SoftOpt}
\max_{a_i, i\in[1,n]} \quad \mu_s, \quad \textmd{subject to} \quad \lambda_p<\mu_p,
\end{equation}
where $\mu_s$ is the secondary user throughput, and $\mu_p$ is the primary user service rate. Next, we characterize $\mu_p$.

Under the assumption stated before that secondary users always have packets to send, the service process of the $q^{th}$ PU can be characterized as
\begin{equation}
Y_q^t=\mathbf{1}\left(A_q^t\bigcap\overline{O_{qd}^t}\bigcap_{l\in\mathcal{M}_s}\left\{\overline{\mathcal{B}\bigcap P_s}\right\}\right),
\end{equation}
where $\mathbf{1}(\cdot)$ denotes the indicator function ($\mathbf{1}(A)=1$ if event $A$ occurs, and $0$ otherwise), $A_q^t$ denotes the event that time slot $t$ is assigned to primary user $q$, $\overline{O_{qd}^t}$ denotes the event that the link between PU $q$ and its destination is not in outage, $\mathcal{B}$ is the event of missed detection, and $P_s$ is the event that a SU gains access to the channel. The probability of the joint event of missed detection and permission to access the channel, denoted $p_s^1$, is given by $Pr\left\{\mathcal{B}\bigcap P_s\right\}=p_s^1=\sum_{i\in[1,n]}p_i^1a_i$,
where $p_i^1$ is the probability that the energy detector's output of the received signal $\left\|y_{ps}\right\|^2$ falls in the $i^{th}$ subinterval when the PU is present. From the received signal model of (\ref{eqn: Rx Signal}), $p_i^1=\exp\left(-\frac{(i-1)\eta}{2n\sigma_1^2}\right)-\exp\left(-\frac{i\eta}{2n\sigma_1^2}\right)$, where $\sigma_1^2$ is the variance of the energy detector's output when the PU is present.

The average PU service rate can now be written as
\begin{equation}\label{eqn: mu_p_noFB}
\mu_p=E\left\{Y_q^t \right\}=\frac{1-P_{pd}^o}{M_p}\left(1-\sum_{i\in[1,n]}p_i^1a_i\right)^{M_s},
\end{equation}
where $E\left\{\cdot\right\}$ is the expectation operator, and $P_{pd}^o$ is the probability that the link between the primary transmitter and the primary receiver is in outage. Next, we move to characterizing $\mu_s$.

For a secondary user to successfully send its packet, the following events have to all take place simultaneously: it has to correctly identify the channel as idle, i.e. no false alarm occurs, it must gain access to the channel, its own link must not be in outage, all other secondary users must either have a false alarm decision or have no access to the channel, and the PU's queue has to be empty. Thus, the service process of the $k^{th}$ SU can be characterized as
\begin{eqnarray}\label{eqn: SU_servEvent}
Y_k^t=\mathbf{1}\Bigg{(}\bigcup_{q\in\mathcal{M}_p}\bigg{[}A_q^t\bigcap\left\{Q_q^t=0\right\}
\bigcap\overline{O_{kd}^t}\bigcap\mathcal{\overline{A}}\bigcap P_s \nonumber \\
 \bigcap_{l\in\mathcal{M}_s\setminus k}\left\{\mathcal{A}\bigcup \overline{P_s}\right\}\bigg{]}\Bigg{)} ,
\end{eqnarray}
where $\mathcal{A}$ is the event of false alarm, and $\{Q_q^t=0\}$ denotes the event that the $q^{th}$ PU queue is empty, which can be determined using Little's theorem \cite{Kleinrock} to be $(1-\lambda_p/\mu_p)$.

The joint event of no false alarm and gaining channel access when the PU is not present can be expressed as $Pr\left\{\mathcal{\overline{A}}\bigcap P_s\right\}=p_s^0=\sum_{i\in[1,n]}p_i^0a_i$,
where $p_i^0=\exp\left(-\frac{(i-1)\eta}{2n\sigma_0^2}\right)-\exp\left(-\frac{i\eta}{2n\sigma_0^2}\right)$, and $\sigma_0^2$ is the variance of the energy detector's output when no PU is present.

Therefore, the average SU service rate is given by
\begin{align}\label{eqn: mu_s_noFB}
\mu_s=&E\left\{Y_k^t \right\} =\Bigg{(}1-\frac{\lambda_pM_p}{\big{(}1-P_{pd}^0\big{)}\big{(}1-\sum_{i\in[1,n]}p_i^1a_i\big{)}^{M_s}}\Bigg{)}\times \nonumber \\
&(1-P_{sd}^o)\Bigg{(}\sum_{i\in[1,n]}p_i^0a_i\Bigg{)}\Bigg{(}1-\sum_{i\in[1,n]}p_i^0a_i\Bigg{)}^{M_s-1}.
\end{align}

Fortunately, the optimization problem (\ref{eqn: SoftOpt}) using (\ref{eqn: mu_p_noFB}) and (\ref{eqn: mu_s_noFB}) was proved to be convex \cite{Soft_Sensing}. Thus, its global optimum can be calculated efficiently via standard techniques \cite{boyd}.

\section{Proposed Feedback- and Soft Sensing-Based Access}\label{FbScheme}

In our proposed scheme, secondary users overhear the primary user feedback, which is assumed to be error-free, by the end of each time slot and leverage it as follows:
\begin{itemize}
\item If an ACK/no feedback is heard, each SU attempts sensing and accessing the channel in the next time slot.
\item Otherwise, if a NACK is heard, all secondary users back-off completely in the next time slot allowing for retransmission of the erroneous primary packet. This, in turn, avoids guaranteed collisions with the PU, increasing the PU service rate and decreasing the primary user's packet delay. Accordingly, the primary queue will be empty with a higher probability which increases the throughput of the secondary network.
\end{itemize}

\begin{figure}
\includegraphics[scale=0.34]{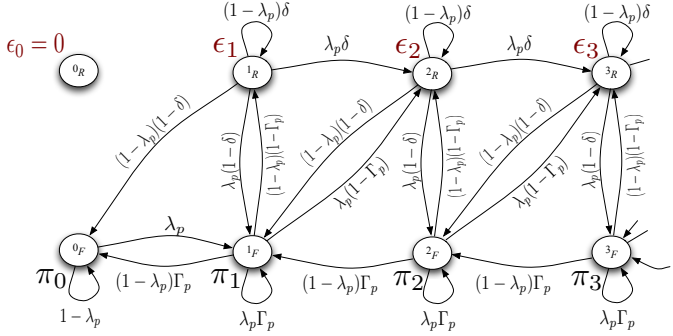}
\centering
\caption{Markov Chain model of the PU queue evolution.}\label{fig: Markov_FB}
\vspace{-0.15in}
\end{figure}

The Markov Chain modeling the PU queue dynamics is presented in Fig. \ref{fig: Markov_FB}. There are two classes of states the PU queue may encounter, the first is $k_F$, denoting the case where the PU has $k$ packets and sending for the first time, where $F$ stands for ``First transmission''. The second is $k_R$, denoting the case where the PU has $k$ packets and re-transmitting, where $R$ stands for ``Retransmission''. States $k_F$ have stationary probability $\pi_k$, and states $k_R$ have stationary probability $\epsilon_k$.

The second class of states, $k_R$, is only reached after the reception of a NACK. In such case, the primary link's outage is the sole cause of packet retransmissions since secondary users refrain completely from transmission upon overhearing a NACK from the primary receiver. If the PU queue is empty, then, clearly the PU cannot be in the retransmission state, therefore, $\epsilon_0=0$.

As shown in Fig. \ref{fig: Markov_FB}, a down transition from $(k+1)_F$ in slot $t$ to $k_F$ in slot $t+1$ occurs when the PU does not receive any packets during slot $t$, which occurs with probability $1-\lambda_p$, and at the same time succeeds in transmission, which occurs with probability $\Gamma_p=\frac{1}{M_p}(1-P_{pd}^o)(1-\sum_{i\in[1,n]}p_i^1a_i)^{M_s}$, i.e. the PU got allocated to this time slot, at the same time its link was not in outage, and all secondary users did not interfere with its transmission. These two events are independent, and hence, their joint probability simply boils down to their product. It is worth noting here that $\Gamma_p$ has the same value as the primary user service rate $\mu_p$ in the baseline no feedback scheme since they both denote the successful primary transmission probability in the same surrounding conditions. A PU will stay in state $k_F$ in time slot $t+1$, for $k\geq1$, if it received a new packet during time slot $t$, which occurs with probability $\lambda_p$, and if it simultaneously succeeded in transmission, which occurs with probability $\Gamma_p$. Again, these two events are independent and hence there joint probability is equal to their multiplication, that is $\lambda_p\Gamma_p$.

On the other hand, an up transition from $k_R$ in slot $t$ to $(k+1)_R$ in slot $t+1$ occurs when the PU receives a packet in slot $t$, which occurs with probability $\lambda_p$, and fails in its transmission, which should now be considered in the absence of secondary users. This occurs with probability $\delta=1-\frac{1}{M_p}(1-P_{pd}^o)$, which means that either the PU did not gain access to the slot, or it gained access, yet, its channel was in outage. It follows that the probability of a PU staying in state $\epsilon_k$ is equal to $(1-\lambda_p)\delta$. The rest of the probabilities can be derived using similar arguments.

Next, we present our system analysis and consider two performance metrics, namely the SU throughput and PU packet delay.

\subsection{Secondary Throughput Analysis}

\subsubsection{Problem Formulation}

In this subsection, we derive an expression for the SU throughput in the proposed feedback-based scheme. The SU service event will be just the same as in (\ref{eqn: SU_servEvent}). Due to the feedback, it is only the value of $Pr\{Q_q^t=0\}$ (which is equivalent in our model to $\pi_0$) that is going to change. The PU Markov chain in Fig. \ref{fig: Markov_FB} can be analyzed, using the global balance equations, in order to get the value of $\pi_0$ which is given by
\begin{equation}\label{eqn: pi_0_FB}
\pi_0=\frac{\chi-\lambda_p}{1-\delta},
\end{equation}
where $\chi=\lambda_p\Gamma_p+(1-\lambda_p)(1-\delta)$ (proof in Appendix). After some algebraic manipulations, we get the following result
\begin{equation}\label{eqn: pi_0_FB_extended}
\pi_0=1-\lambda_p\left[\left(1+\frac{M_p}{1-P_{pd}^o}\right)-\left(1-\sum_{i\in[1,n]}p_i^1a_i\right)^{M_s}\right].
\end{equation}
It is worth noting that for an irreducible and aperiodic Markov chain, the queue is stable if there exists a non-zero value for the probability of the queue being empty \cite{Szpan94} . This condition is equivalent in our model to having $\pi_0>0$, which leads to $\lambda_p<\chi$ (which is the same condition of $\psi<1$ stated in (\ref{eqn: term_B}) in Appendix).

In order to gain more insight into the difference in the SU throughput between our proposed scheme and the no feedback one, we compute $\triangle\pi_0$, denoting the difference between $\pi_0$ in the feedback-based scheme (equation (\ref{eqn: pi_0_FB})), and $\pi_0$ in the no feedback one that turns out to be equal to $1-\lambda_p/\mu_p$ (directly from Little's law), which is equivalent to $1-\lambda_p/\Gamma_p$, since $\mu_p$ and $\Gamma_p$ have the same value. Therefore, and after some algebraic manipulations, the following result can be reached
\begin{equation}
\triangle\pi_0=\frac{\lambda_p(1-\Gamma_p)}{\Gamma_p}\left[1-\left(1-\sum_{i\in[1,n]}p_i^1a_i\right)^{M_s}\right],
\end{equation}
which is always positive for $0\leq a_i\leq1$. Therefore, the SU throughput of the proposed feedback-based scheme is always larger than that of the no feedback one for the same set of SU access probabilities. One can expect that finding access probabilities that maximize the throughput for the feedback-based scheme should give even higher SU throughput.

We can now write the formula of the SU throughput as
\begin{align}\label{eqn: mu_s_FB}
\mu_s=&\Bigg{(}1-\lambda_p\Bigg{[}\bigg{(}1+\frac{M_p}{1-P_{pd}^o}\bigg{)}-\bigg{(}1-\sum_{i\in[1,n]}p_i^1a_i\bigg{)}^{M_s}\Bigg{]}\Bigg{)}\times
\nonumber \\
&(1-P_{sd}^0)\Bigg{(}\sum_{i\in[1,n]}p_i^0a_i\Bigg{)}\Bigg{(}1-\sum_{i\in[1,n]}p_i^0a_i\Bigg{)}^{M_s-1}.
\end{align}
Therefore, the optimization problem is given by
\begin{equation}\label{eqn: FBOpt}
\max_{a_i, i\in[1,n]} \quad \mu_s, \quad \textmd{subject to} \quad \lambda_p<\chi.
\end{equation}
Unfortunately, the optimization problem in (\ref{eqn: FBOpt}) is non-convex as we show later. Nevertheless, It can still be solved efficiently by exploiting its structure as discussed next.

\subsubsection{General Optimization Approach}

Consider the following maximization problem
\begin{equation}\label{eqn: max_toy_exmpl}
\max_{x} \quad f_1(x) + f_2(x), \quad s.t. \quad 0\leq x \leq 1
\end{equation}
where $f_1(x)$ is a concave function in $x$, while $f_2(x)$ is non concave. Now let us assume without loss of generality, that the function $f_2(x)$ is bounded from below and from above by $f_2^{min}$ and $f_2^{max}$ respectively for any given value of $x$ in the feasible region.

Now consider the following algorithm
\begin{align*}
&\quad \quad \texttt{OPTIMIZATION ALGORITHM} \nonumber \\
&\texttt{sum}=-\infty \nonumber \\
&\textbf{LOOP}: \quad \tau=f_2^{min}:\nu:f_2^{max} \nonumber \\
&\quad \max_{x} \quad f_1(x) \nonumber \\
&\quad s.t. \quad 0\leq x \leq 1 \nonumber \\
&\quad \quad \quad f_2(x)\geq\tau, \nonumber \\
&\quad \texttt{dummy}=f_1(x^*)+f_2(x^*) \nonumber \\
&\quad \textbf{if} \quad \{\texttt{dummy} \geq \texttt{sum}\} \nonumber \\
&\quad \quad \texttt{sum}=\texttt{dummy} \nonumber \\
&\quad \textbf{end if} \nonumber \\
&\textbf{end LOOP}
\end{align*}
where $\nu$ is the step size, and $x^*$ is the value of $x$ that maximizes $f_1(x)$ while satisfying the constraints. As $\nu\rightarrow0$, the algorithm introduced here gives the same solution for (\ref{eqn: max_toy_exmpl}). The proof of this is straightforward. In the search for the optimal value of the variable $x$, it is made sure that $f_1$ is maximized and, simultaneously, $f_2$ has a value larger than or equal to $\tau$, which is an iteration term that takes on the possible values of the bounded function $f_2$. Once the problem is solved, the value $f_1+f_2$ is computed and compared to the largest saved value. If it is larger, it is then saved as the new largest value. If not, the algorithm continues to the next iteration. Eventually, the largest value is reached. This reformulation is relatively efficient to solve if, for each $\tau$, the optimization problem inside the loop is convex. This requires that the function $f_1$ is concave and the inequality constraint $f_2>\tau$ can be cast in the form of a concave function greater than a constant \cite{boyd}. We next show that (\ref{eqn: FBOpt}) can be solved using the aforementioned algorithm.

\subsubsection{Solution Approach}

First, we take the logarithm of the expression of $\mu_s$ in (\ref{eqn: mu_s_FB}) before the maximization. This yields an equivalent problem that has the same solution since the $\log(.)$ is a monotonic function. The expression now becomes
\begin{align}\label{eqn: mu_s_FB_log}
&\log(\mu_s)= \nonumber \\
&\log\underbrace{\Bigg{(}1-\lambda_p\Bigg{[}\bigg{(}1+\frac{M_p}{1-P_{pd}^o}\bigg{)}-\bigg{(}1-\sum_{i\in[1,n]}p_i^1a_i\bigg{)}^{M_s}\Bigg{]}\Bigg{)}}_{\pi_0} \nonumber \\
&+\log(1-P_{sd}^0)+\log\Bigg{(}\sum_{i\in[1,n]}p_i^0a_i\Bigg{)} \nonumber \\ &+(M_s-1)\log\Bigg{(}1-\sum_{i\in[1,n]}p_i^0a_i\Bigg{)}.
\end{align}
The last two terms in (\ref{eqn: mu_s_FB_log}) are the logarithm of an affine function in $a_i$ and hence are concave in $a_i$ \cite{boyd}.

The first term, however, is the logarithm of a convex function in $a_i$ (proof is omitted due to space limitations), which causes a problem since the logarithm of a convex function cannot be proven to be concave \cite{boyd}. But since $\pi_0$ (the term inside the logarithm) is bounded between zero and one, we can use the \textit{Optimization Algorithm} presented above to solve this optimization problem as follows. First, we divide the term $\log(\mu_s)$ into two parts, the first consists of the sum of all the concave terms, and the second consists of the non-concave term $\log(\tau)$, where $\tau=\pi_0$. The two parts map into $f_1(x)$ and $f_2(x)$ in the \textit{Optimization Algorithm} presented above, respectively. The second step is to apply the algorithm as follows
\begin{align}
&\texttt{sum}=-\infty \nonumber\\
&\textbf{LOOP}: \quad \tau=0:\nu:1 \nonumber\\
&\max_{a_i, i\in[1,n]} \log\bigg{(}\sum_{i\in[1,n]}p_i^0a_i\bigg{)}+(M_s-1)\log\bigg{(}1-\sum_{i\in[1,n]}p_i^0a_i\bigg{)} \label{eqn: FBOpt_after_mod}\\
&\quad s.t. \quad 0\leq a_i \leq 1, \quad \forall i\in[1,n] \nonumber\\
&\quad \quad \quad \log\bigg{(}1-\sum_{i\in[1,n]}p_i^1a_i\bigg{)}\geq\frac{1}{M_s}\log\big{[}f(\tau)\big{]}^+\label{eqn: newly_added_constraint}\\
&\quad \texttt{dummy}= \mu_s(a_i^*)\nonumber\\
&\quad \textbf{if} \quad \{\texttt{dummy} \geq \texttt{sum}\} \nonumber\\
&\quad \quad \texttt{sum}=\texttt{dummy} \nonumber\\
&\quad \textbf{end if} \nonumber\\
&\textbf{end LOOP}, \nonumber
\end{align}
where $f(\tau)=\frac{\tau}{\lambda_p}+\frac{M_p}{1-P_{pd}^o}-\frac{1-\lambda_p}{\lambda_p}$, and $\left[f(\tau)\right]^+$ denotes $\max(0,f(\tau))$. Accordingly, (\ref{eqn: FBOpt_after_mod}) is now concave and can be solved using standard convex optimization tools. It must be noted that the stability condition in (\ref{eqn: FBOpt}) can be rewritten as
\begin{equation}
\log\bigg{(}1-\sum_{i\in[1,n]}p_i^1a_i\bigg{)}\geq\frac{1}{M_s}\log\bigg{(}\bigg{[}\frac{M_p}{1-P_{pd}^o}-\frac{1-\lambda_p}{\lambda_p}\bigg{]}^+\bigg{)},
\end{equation}
which is subsumed by the newly added constraint (\ref{eqn: newly_added_constraint}) as it corresponds to $\tau=0$ and $\tau$ is non-negative. Therefore, we have managed to overcome the problem of the non-convexity of the optimization problem in (\ref{eqn: FBOpt}) via a simple algorithm which requires an exhaustive search over only one bounded parameter $\tau$.

\subsection{Primary Delay Analysis}

In this subsection, we only present final expressions for the average PU packet delay. Proofs of these are omitted due to space limitations. For the no feedback soft sensing scheme presented in Section \ref{SoftSense}, one can easily show that
\begin{equation}\label{eqn: Delay_noFB}
D_p=\frac{1-\lambda_p}{\mu_p-\lambda_p}.
\end{equation}
While for the proposed feedback-based scheme, the delay is given by
\begin{equation}\label{eqn: Delay_FB}
D_p=\frac{(\Gamma_p-\chi)(\chi-\lambda_p)^2+(1-\lambda_p)^2(1-\Gamma_p)\chi}{(1-\lambda_p)(1-\chi)(1-\delta)(\chi-\lambda_p)}.
\end{equation}

\section{Performance Results}\label{Simulations}

In this section, we compare the performance of our proposed feedback-based scheme with two other schemes, namely the conventional (non-feedback-based) soft sensing scheme, and the Neyman-Pearson hard decision scheme. We consider a system of $M_p=4$ primary users and $M_s=2$ secondary users. The distance between the primary transmitters and receivers is set to 100 m, the distance between the secondary transmitters and receivers is also set to 100 m, and the distance between any primary user and any secondary user is set to 150 m. The SNR threshold $\zeta$ is 10 dB, the transmit power is 100 mW, the path loss exponent $\gamma=3.7$, and $N_0=10^{-11}$ W/Hz. The region below the energy threshold $\eta$ is divided into $n=4$ regions each having a different access probability.

\begin{figure}
\includegraphics[scale=0.43]{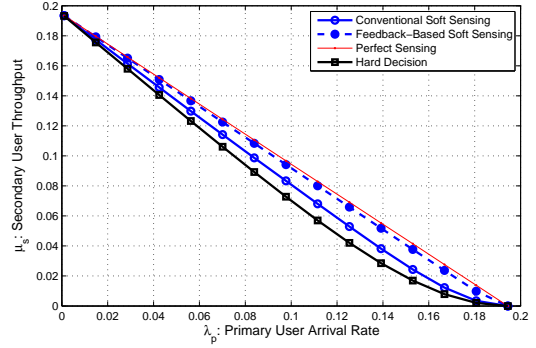}
\centering
\caption{Comparison between SU throughput of different schemes in a system of 4 PUs and 2 SUs.}\label{fig: SUthroughput}
\vspace{-0.15in}
\end{figure}

In Fig. \ref{fig: SUthroughput}, the SU throughput is plotted against the PU arrival rate. Different schemes are compared with respect to the upper bound acquired by perfect sensing; a scheme that can be considered genie-aided, where the SU perfectly knows when the PU is idle. We can see that our proposed feedback-based scheme, when applied jointly with soft sensing, not only outperforms the conventional soft sensing one but also approaches the upper bound almost with equality in some regions. Also the Neyman-Pearson (N-P) hard decision sensing scheme is plotted for completeness.

\begin{figure}
\includegraphics[scale=0.43]{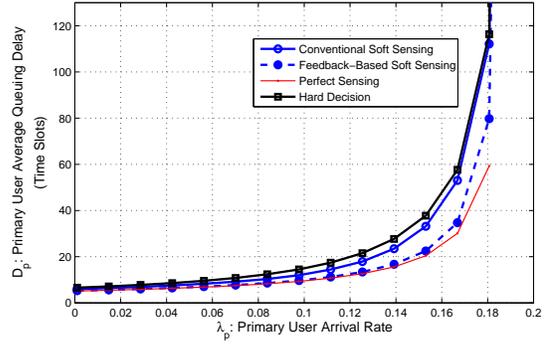}
\centering
\caption{Comparison between PU packet delay of different schemes in a system of 4 PUs and 2 SUs.}\label{fig: PUDelay}
\vspace{-0.15in}
\end{figure}

In Fig. \ref{fig: PUDelay}, the average PU queuing delay is plotted against the PU arrival rate. We can see that our proposed feedback-based scheme when applied jointly with soft sensing also outperforms the conventional soft sensing and the hard decision sensing ones in all regions.

\begin{figure}
\includegraphics[scale=0.43]{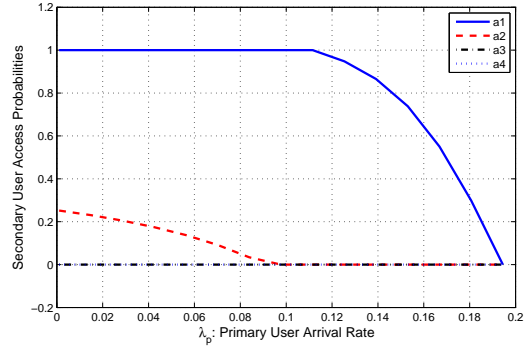}
\centering
\caption{SU access probabilities in a conventional soft sensing scheme in a system of 4 PUs and 2 SUs.}\label{fig: AccessProbNoFB}
\vspace{-0.15in}
\end{figure}

\begin{figure}
\includegraphics[scale=0.43]{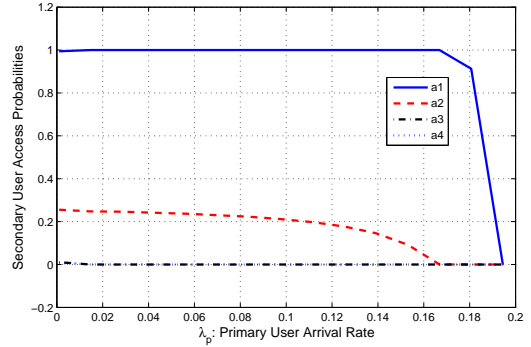}
\centering
\caption{SU access probabilities in a feedback-based soft sensing scheme in a system of 4 PUs and 2 SUs.}\label{fig: AccessProbFB}
\vspace{-0.15in}
\end{figure}

Access probabilities for both the conventional soft sensing scheme and for the feedback-based one are plotted against the PU arrival rate in Fig. \ref{fig: AccessProbNoFB} and Fig. \ref{fig: AccessProbFB}, respectively. From the figures, we can see that in both cases the two access probabilities closer to the decision threshold $a_3$ and $a_4$ are equal to zero for any given arrival rate. However, $a_1$ and $a_2$ have higher values in the feedback-based scheme than their counterparts in the conventional soft sensing one at relatively high arrival rates. This is attributed to the proper use of the primary feedback information by the secondary users, which makes them avoid sure collisions, and thus enables them to access the channel more aggressively without affecting the PU's stability.

\begin{figure}
\includegraphics[scale=0.4497]{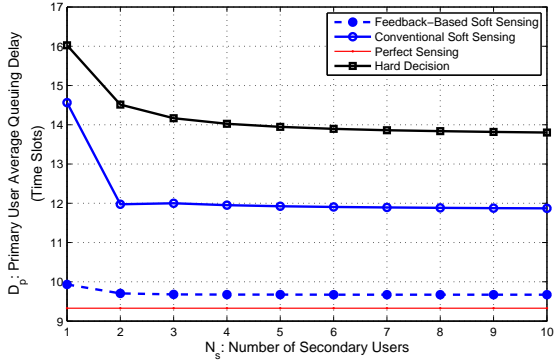}
\centering
\caption{PU packet delay vs. number of secondary users at $\lambda_p\simeq0.1$ and $M_p=4$ PUs.}\label{fig: PUdelay_multiSU}
\vspace{-0.15in}
\end{figure}

In order to gain more insights into how our proposed scheme performs with different number of secondary users, we provide scalability results. For an arrival rate of $\lambda_p\simeq0.1$ packets per time slot, a plot of the PU packet delay for different schemes against the number of secondary users is presented in Fig. \ref{fig: PUdelay_multiSU}. We can see that our proposed feedback-based scheme is the nearest to the lower bound at any given number of secondary users. We also notice that the PU packet delay curve converges to a certain level. This is due to the fact that the access probabilities change inversely proportional to the number of secondary users in order to guarantee the stability of the primary users' queues. This opposite change also causes the PU service rate to converge to a certain level, thereby causing the delay to be constant.

\begin{figure}
\includegraphics[scale=0.4497]{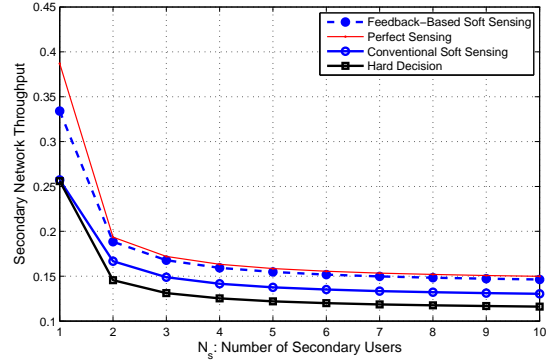}
\centering
\caption{Secondary network throughput vs. number of secondary users at $\lambda_p\simeq0.1$ and $M_p=4$ PUs.}\label{fig: NWthroughput_multiSU}
\vspace{-0.15in}
\end{figure}

Fig. \ref{fig: NWthroughput_multiSU} presents a result pertaining to the secondary network throughput $(M_s*\mu_s)$ against the number of secondary users, for $\lambda_p\simeq0.1$ packets per time slot too. The network throughput seems to be slowly decreasing with the increase of secondary users, however, our proposed feedback-based scheme outperforms both the conventional soft sensing and hard decision schemes at every given number of secondary users and closely approaches the optimal perfect sensing scheme.

\section{Conclusions}\label{Concl}

In this paper, we examined a cognitive spectrum access scheme in which secondary users exploit the primary feedback information. We considered a secondary network employing a random access scheme in which secondary users access the channel by certain access probabilities that are function of the sensing metric. We studied the cognitive radio network from a queuing theory point of view. Access probabilities are determined by solving an optimization problem subject to a constraint on the primary user queue stability. In setting our problem, we assumed that secondary users can eavesdrop on the primary link's feedback; secondary users back-off completely from accessing the channel upon hearing a NACK, and attempt accessing if an ACK/no feedback is overheard. This has led to significant results in both secondary user throughput and primary user packet delay. Our proposed scheme has outperformed both the soft sensing and the conventional hard decision sensing schemes and produced very close results to the \textit{optimal} genie-aided scheme in which secondary users have perfect knowledge of the activity of the primary user and act upon it.

\appendix


Referring to the Markov chain in Fig. \ref{fig: Markov_FB}, we can write the global balance equation around state $0_F$ as follows
\begin{equation}\label{eqn: anls_pi0}
\pi_0\lambda_p=\pi_1\bar{\lambda}_p\Gamma_p+\epsilon_1\bar{\lambda}_p\bar{\delta},
\end{equation}
where the notation $\bar{x}=1-x$. Writing the balance equation around state $1_R$ we get
\begin{equation*}
\epsilon_1(1-\delta\bar{\lambda}_p)=\pi_1\bar{\lambda}_p\bar{\Gamma}_p,
\end{equation*}
therefore, we have
\begin{equation}\label{eqn: anls_eps1}
\pi_1=\epsilon_1\frac{1-\delta\bar{\lambda}_p}{\bar{\lambda}_p\bar{\Gamma}_p}.
\end{equation}
Substituting by (\ref{eqn: anls_eps1}) in (\ref{eqn: anls_pi0}), we get
\begin{equation}\label{eqn: anls_eps1_pi0}
\epsilon_1=\frac{\lambda_p\bar{\Gamma}_p}{\chi}\pi_0,
\end{equation}
where $\chi=\lambda_p\Gamma_p+\bar{\lambda}_p\bar{\delta}$. Now using (\ref{eqn: anls_eps1_pi0}) in (\ref{eqn: anls_eps1}) yields
\begin{equation}\label{eqn: anls_pi1_pi0}
\pi_1=\frac{\lambda_p(1-\delta\bar{\lambda}_p)}{\bar{\lambda}_p\chi}\pi_0.
\end{equation}
Writing the balance equation around state $1_F$, we have
\begin{equation*}
\pi_1(1-\lambda_p\Gamma_p)=\pi_0\lambda_p+\epsilon_1\lambda_p\bar{\delta}+\pi_2\bar{\lambda}_p\Gamma_p+\epsilon_2\bar{\lambda}_p\bar{\delta}.
\end{equation*}
Using (\ref{eqn: anls_pi0}) to substitute for the term $\pi_0\lambda_p$, we get
\begin{equation}\label{eqn: anls_pi1}
\pi_1\bar{\Gamma}_p=\epsilon_1\bar{\delta}+\pi_2\bar{\lambda}_p\Gamma_p+\epsilon_2\bar{\lambda}_p\bar{\delta}.
\end{equation}
Using (\ref{eqn: anls_eps1_pi0}) and (\ref{eqn: anls_pi1_pi0}) into (\ref{eqn: anls_pi1}), we now have
\begin{equation}\label{eqn: anls_pi2_eps2_pi0}
\pi_2\bar{\lambda}_p\Gamma_p+\epsilon_2\bar{\lambda}_p\bar{\delta}=\frac{\lambda_p^2\bar{\Gamma}_p}{\bar{\lambda}_p\chi}\pi_0.
\end{equation}
Writing the balance equation around state $2_R$, we get
\begin{equation*}
\epsilon_2(1-\delta\bar{\lambda}_p)=\epsilon_1\lambda_p\delta+\pi_1\lambda_p\bar{\Gamma}_p+\pi_2\bar{\lambda}_p\bar{\Gamma}_p,
\end{equation*}
but since from (\ref{eqn: anls_eps1_pi0}) and (\ref{eqn: anls_pi1_pi0}) we have
\begin{equation*}
\epsilon_1\lambda_p\delta+\pi_1\lambda_p\bar{\Gamma}_p=\frac{\lambda_p^2\bar{\Gamma}_p}{\bar{\lambda}_p\chi}\pi_0,
\end{equation*}
therefore
\begin{equation}\label{eqn: anls_eps2}
\epsilon_2(1-\delta\bar{\lambda}_p)-\pi_2\bar{\lambda}_p\bar{\Gamma}_p=\frac{\lambda_p^2\bar{\Gamma}_p}{\bar{\lambda}_p\chi}\pi_0.
\end{equation}
From (\ref{eqn: anls_pi2_eps2_pi0}) and (\ref{eqn: anls_eps2}) we can get the following
\begin{equation}\label{eqn: anls_pi2_eps2}
\epsilon_2=\frac{\bar{\lambda}_p}{\lambda_p}\pi_2.
\end{equation}
Therefore, using (\ref{eqn: anls_pi2_eps2}) in (\ref{eqn: anls_pi2_eps2_pi0}) we get
\begin{equation}\label{eqn: anls_k2_pi0}
\epsilon_2=\big{(}\frac{\lambda_p\bar{\chi}}{\bar{\lambda}_p\chi}\big{)}^2.\frac{\bar{\lambda}_p\bar{\Gamma}_p}{\bar{\chi}^2}\pi_0, \quad \textmd{and} \quad \pi_2=\big{(}\frac{\lambda_p\bar{\chi}}{\bar{\lambda}_p\chi}\big{)}^2.\frac{\lambda_p\bar{\Gamma}_p}{\bar{\chi}^2}\pi_0.
\end{equation}

From the symmetry of the upcoming states in the Markov chain, one can expect that equation (\ref{eqn: anls_pi2_eps2}) can be generalized for any $\epsilon_k$ and $\pi_k$ with $k\geq2$, since all the upcoming balance equations will give the same result. Also this applies for the results in (\ref{eqn: anls_k2_pi0}). Verification of this is straight forward but it is omitted due to space limits.

Therefore, we can now write the following results:
\begin{itemize}
\item $\epsilon_0=0.$
\item $\epsilon_1=\frac{\lambda_p\bar{\Gamma}_p}{\chi}\pi_0.$
\item $\pi_1=\frac{\lambda_p(1-\delta\bar{\lambda}_p)}{\bar{\lambda}_p\chi}\pi_0.$
\end{itemize}
And for $k\geq2$ we have:
\begin{itemize}
\item $\epsilon_k=\big{(}\frac{\lambda_p\bar{\chi}}{\bar{\lambda}_p\chi}\big{)}^k.\frac{\bar{\lambda}_p\bar{\Gamma}_p}{\bar{\chi}^2}\pi_0.$
\item $\pi_k=\frac{\lambda_p}{\bar{\lambda}_p}\epsilon_k.$
\end{itemize}

We can now use the normalization condition, $\sum_{k=0}^{\infty}(\pi_k+\epsilon_k)=1$, to get the value of $\pi_0$. First, we will divide the summation as follows
\begin{equation}\label{eqn: normCond}
\sum_{k=0}^{\infty}(\pi_k+\epsilon_k)=\pi_0+\underbrace{(\pi_1+\epsilon_1)}_A+\underbrace{\sum_{k=2}^{\infty}(\pi_k+\epsilon_k)}_B=1.
\end{equation}
Simplifying the term $B$: since, for $k\geq2$, we have
\begin{equation*}
\pi_k+\epsilon_k=\psi^k\frac{\bar{\Gamma}_p}{\bar{\chi}^2}\pi_0, \quad \textmd{where} \quad \psi=\frac{\lambda_p\bar{\chi}}{\bar{\lambda}_p\chi}.
\end{equation*}
Hence,
\begin{equation}\label{eqn: term_B}
B=\frac{\bar{\Gamma}_p\pi_0}{\bar{\chi}^2}\sum_{k=2}^{\infty}\psi^k=\bigg{(}\frac{\lambda_p\bar{\Gamma}_p}{\bar{\lambda}_p\chi}\bigg{)}\bigg{(}\frac{\lambda_p}{\chi-\lambda_p}\bigg{)}\pi_0.
\end{equation}
The last summation converges only if $\psi<1$, that is equivalent to $\lambda_p<\chi$. This is actually the stability condition for the PU queue. After some manipulations, the term $A$ can be written as:
\begin{equation}\label{eqn: term_A}
A=\bigg{(}\frac{\lambda_p\bar{\Gamma}_p}{\bar{\lambda}_p\chi}\bigg{)}\bigg{(}\frac{\chi+\bar{\Gamma}_p}{\bar{\Gamma}_p}\bigg{)}\pi_0.
\end{equation}
From (\ref{eqn: term_B}) and (\ref{eqn: term_A}), and after some involved manipulations, the final result becomes
\begin{equation}\label{eqn: term_AplusB}
A+B=\frac{\lambda_p(\bar{\Gamma}_p+\bar{\delta})}{\chi-\lambda_p}\pi_0.
\end{equation}
Using this final result of (\ref{eqn: term_AplusB}) in (\ref{eqn: normCond}), we can write the value of $\pi_0$ as
\begin{equation}\label{eqn: anls_pi0_FB}
\pi_0=\frac{\chi-\lambda_p}{\bar{\delta}},
\end{equation}
which can be checked to satisfy the balance equation given in (\ref{eqn: anls_pi0}).

\bibliographystyle{IEEEbib}
\bibliography{MyLib}

\end{document}